# A Pattern Theorem for Lattice Clusters

Neal Madras*
Department of Mathematics and Statistics
York University
4700 Keele St.
Toronto, Ontario, M3J 1P3 CANADA
E-mail: madras@nexus.yorku.ca

July 6, 1999

**Abstract**

We consider general classes of lattice clusters, including various kinds of animals and trees on different lattices. We prove that if a given local configuration ("pattern") of sites and bonds can occur in large clusters, then it occurs at least $cN$ times in most clusters of size $n$, for some constant $c > 0$. An analogous theorem for self-avoiding walks was proven in 1963 by Kesten. We use the pattern theorem to prove the convergence of $\lim_{n \to \infty} a_{n+1}/a_n$ (where $a_n$ is the number of clusters of size $n$, up to translation). The results also apply to weighted sums, and in particular we can take $a_n$ to be the probability that the percolation cluster containing the origin consists of exactly $n$ sites. Another consequence is strict inequality of connective constants for sublattices and for certain subclasses of clusters.

## 1 Introduction

Let $L$ be a periodic lattice in $d$-dimensional Euclidean space $\mathbf{R}^d$ ($d \geq 2$). We shall give a precise definition of "periodic lattice" in Section **??** but for now we'll just think of some common examples: the hypercubic lattice $\mathbf{Z}^d$ of integral points with nearest-neighbour bonds, the triangular and hexagonal lattices in two dimensions, and the face-centered and body-centered cubic lattices in three dimensions. We shall think of a lattice as a graph with infinitely many sites and bonds embedded in $\mathbf{R}^d$. The above examples are all undirected graphs, and we shall restrict our attention to undirected graphs for most of the paper. Our main results do extend to directed lattices as well, but there are some technical differences; see Section 3.6 for a discussion.

*This work was supported in part by the Natural Sciences and Engineering Research Council of Canada. The author was visiting the Fields Institute for Research in Mathematical Sciences, in Toronto, while writing of the paper.

This paper concerns "clusters" of the lattice, which is a generic term we use to denote one of several possible families of <u>finite</u> subgraphs of $L$, including:

- <u>Bond animals</u>: A bond animal is simply a finite connected subgraph of $L$.

- <u>Site animals</u>: A site animal is a finite connected subgraph $G$ of $L$ with the property that if $b$ is bond of $L$ that has both endpoints in $G$, then $b$ must be in $G$. Thus a site animal is determined by its set of sites, in contrast to a bond animal.

- <u>Bond trees</u>: A bond tree is a bond animal with no cycles.

- <u>Directed clusters</u>: A directed cluster is usually described as a cluster on a directed lattice with the property that each site of the cluster can be reached by a directed path starting from a fixed site of the cluster (the "root"). We will take a view that is slightly less general, but still broad enough to include most standard examples.

Unfortunately, the methods of this paper do not seem easy to extend to site trees (site animals with no cycles), nor to random surfaces (Vanderzande 1998, Chapter 11).

We are interested in the number of clusters of a given size in a lattice. The word "size" could refer to the number of sites or the number of bonds in the cluster, or perhaps to some other quanitity (see Section 2.2). Since $L$ is infinite and periodic, it makes sense to enumerate clusters up to translation. It has been proven that the number of clusters (up to translation) typically grows exponentially in $n$. Formally, let $C_n$ be the set of all clusters of $L$ of size $n$, and let $C_n^*$ be a subset of $C_n$ that contains exactly one translation of each cluster in $C_n$. In the case of $L = \mathbf{Z}^d$, one can let $C_n^*$ be the set of all clusters in $C_n$ whose lexicographically smallest site is the origin. (For more general lattices, see Section 2.2.) For the kinds of clusters mentioned above (bond or site animals or bond trees, counting by bonds or by sites), one can use concatenation and subadditivity arguments (e.g., Klarner (1967), Klein (1981)) to prove the existence of a *growth constant* $\lambda$, with $1 < \lambda < \infty$, such that

$$\lim_{n \to \infty} |C_n^*|^{1/n} = \lambda. \tag{1}$$

The growth constant $\lambda$ depends on the lattice $L$ as well as on which kind of cluster we are counting. Equation (??) is useful, but it is often too crude to elucidate the physically interesting properties of the model. There is much theoretical and numerical evidence for the belief that there is a *critical exponent* $\theta$ such that

$$|C_n^*| \sim \text{Const}.n^{-\theta} \lambda^n. \tag{2}$$

The physical importance of the critical exponent can be summarized by the belief that, for the kinds of clusters mentioned above, $\theta$ depends only on the dimension of the lattice $L$. For example, bond trees on $\mathbf{Z}^2$ and site animals on the two-dimensional triangular lattice have the same value of $\theta$ (which incidentally is believed to be 1), even though their growth constants are different. This conjecture essentially is a result of the assertion that trees and animals are in the same *universality class*. Hara and Slade (1992) gave a rigorous proof of this for sufficiently high dimensions, but their methods cannot work below eight dimensions. For more on this critical exponent, including what is rigorously known and what is conjectured, see Madras (1995).

One of the main goals of the present paper is to prove a result on the asymptotics of the number of animals that is better than Equation (??) but weaker than (??), namely:

$$\lim_{n \to \infty} \frac{|C^*_{n+1}|}{|C^*_n|} = \lambda. \tag{3}$$

Besides simple enumeration, it is often important to consider sums of weights associated with the clusters. One example of this arises in the *percolation* model (Grimmett, 1989); another arises in collapse models for polymers (Vanderzande 1998, Chapter 8). To describe the collapse model, let $G$ be a bond animal, and let $b$ be a bond of $L$ that is *not in $G$*. We say that $b$ is a *monomer contact* of $G$ if both endpoints of $b$ are sites of $G$, and we say that $b$ is a *solvent contact* of $G$ if exactly one endpoint of $b$ is in $G$. Let mono($G$) and solv($G$) respectively denote the number of monomer contacts and the number of solvent contacts in $G$. Fix two positive parameters $z_m$ and $z_s$, and let the *weight* of $G$ be

$$wt(G) = z_m^{\mathrm{mono}(G)} z_s^{\mathrm{solv}(G)}. \tag{4}$$

Then for each $n$ let

$$\mathcal{G}_n = \sum_{G \in C^*_n} wt(G). \tag{5}$$

Think of an animal $G$ as representing an isolated branched polymer in a solution. Each site of $G$ represent a monomer, and each bond of $G$ is a chemical bond in the polymer. Each solvent contact represents an interaction of a monomer with ions of the solvent, and each monomer contact represents an interaction between two monomers that are close but not connected by a chemical bond. If the the interaction is energetically favourable, then the corresponding $z$ parameter is greater than 1. Viewed as a function of $z_m$ and $z_s$, $\mathcal{G}_n$ is called the *partition function*, and the limit

$$\lim_{n \to \infty} \frac{1}{n} \log \mathcal{G}_n = F(z_m, z_s) \tag{6}$$

is called the *limiting reduced free energy*. This limit is known to exist (Madras et al. (1990), Janse van Rensburg and Madras (1997)). A good understanding of the limiting function $F$, in particular of its differentiability and analyticity properties, gives information about phase transitions in the polymer model.

One can view percolation as a special case of the polymer contact model. Let $p$ be a number between 0 and 1. In standard bond percolation, each bond of the lattice is independently "open" or "closed" with probability $p$ and $1-p$ respectively. Let $\mathcal{C}$ denote the connected component containing the origin in the subgraph of $L$ consisting of all sites and only the open bonds. Then $\mathcal{C}$ is a random subgraph of $L$ (possibly infinite). If $G$ is a fixed bond animal containing the origin, having $n$ sites and $\mathrm{bond}(G)$ bonds, then the probability that $\mathcal{C}$ equals $G$ is

$$\Pr{}_p\{\mathcal{C} = G\} = p^{\mathrm{bond}(G)}(1-p)^{\mathrm{mono}(G)+\mathrm{solv}(G)}, \tag{7}$$

and the probability that $\mathcal{C}$ equals some translation of $G$ is simply $n$ times this quantity. On $\mathbf{Z}^d$ we know that $2dn = 2\,\mathrm{bond}(G) + 2\,\mathrm{mono}(G) + \mathrm{solv}(G)$, so summing over all $n$-site bond animals $G$ containing the origin gives

$$\begin{aligned}
\Pr{}_p&\{\mathcal{C} \text{ has exactly } n \text{ sites}\} \\
&= n \sum_{G \in C_n^*} p^{\mathrm{bond}(G)}(1-p)^{\mathrm{mono}(G)+\mathrm{solv}(G)} \\
&= n \sum_{G \in C_n^*} p^{dn-\mathrm{mono}(G)-\mathrm{solv}(G)/2}(1-p)^{\mathrm{mono}(G)+\mathrm{solv}(G)} \\
&= np^{dn} \sum_{G \in C_n^*} \left(\frac{1-p}{p}\right)^{\mathrm{mono}(G)} \left(\frac{1-p}{\sqrt{p}}\right)^{\mathrm{solv}(G)} \\
&= np^{dn}\mathcal{G}_n,
\end{aligned} \tag{8}$$

where in the last line $\mathcal{G}_n$ is defined as in Equations (??) and (??) using

$$z_m = \frac{1-p}{p} \quad \text{and} \quad z_s = \frac{1-p}{\sqrt{p}}. \tag{9}$$

The analogue of result (??) for the collapse model is

$$\lim_{n \to \infty} \frac{\mathcal{G}_{n+1}}{\mathcal{G}_n} = \exp(F(z_m, z_s)), \tag{10}$$

while for the special case of percolation we have the following result.

**Theorem 1.1** *Consider bond percolation on $\mathbf{Z}^d$ with parameter $p \in (0,1)$. Let $P_p(n)$ be the probability that the open cluster containing the origin has exactly $n$ sites. Then*

$$\lim_{n \to \infty} \frac{P_p(n+1)}{P_p(n)}$$

*exists and equals $\lim_{n \to \infty}[P_p(n)]^{1/n}$.*

We shall see that Theorem ?? also holds for many other lattices, as well as for site percolation and directed percolation (see Theorem ??, Corollary ??, and the end of Section 3.6).

These results are deduced from a different kind of result, known as a "pattern theorem". Roughly speaking, a pattern theorem says that if a certain local configuration of bonds and sites can occur in the middle of a large cluster, then this configuration must occur many times on almost all large clusters. Here, "almost all" means "except for an exponentially small fraction". For example, consider bond animals in $\mathbf{Z}^2$. One local configuration ("pattern") of interest could be a two-by-two lattice square with all eight perimeter bonds present, but the center site (and its four incident bonds) absent. The pattern theorem tells us that there exist positive numbers $\alpha$ and $\epsilon$ such that for all sufficiently large $n$, the fraction of $n$-site bond animals having fewer than $\epsilon n$ occurrences of this particular pattern is at most $e^{-\alpha n}$. More generally, we prove a weighted version of the pattern theorem, where clusters are counted according to their weights. For example, for bond percolation, the pattern theorem becomes a statement that certain probabilities, conditioned on the origin's cluster containing exactly $n$ sites, are exponentially small. Patterns and the pattern theorem are described more fully in Section ??.

The prototype for the results of this paper appeared in Kesten (1963), where ratio limit theorems and pattern theorems were proven for self-avoiding walks. (An $n$-step self-avoiding walk on $L$ is a sequence $\omega_0, \ldots, \omega_n$ of distinct sites of $L$ such that consecutive sites are joined by bonds of $L$.) Let $c_n$ be the number of $n$-step walks starting at the origin. Hammersley and Morton (1954) proved the existence of the limit $\mu = \lim_{n \to \infty} c_n^{1/n}$ on fairly general lattices. For self-avoiding walks on $\mathbf{Z}^d$, Kesten (1963) proved a pattern theorem and deduced that $\lim_{n \to \infty} c_{n+2}/c_n = \mu^2$. The result $\lim_{n \to \infty} c_{n+1}/c_n = \mu$ remains unproven in $\mathbf{Z}^d$, although Kesten's method can be used to prove this stronger result for non-bipartite lattices (see Chapter 7 of Madras and Slade (1993) for further discussion). Hammersley has found a different proof for the pattern theorem for self-avoiding walks (the proof is unpublished, but the method was used in Janse van Rensburg *et al.* (1996) for a similar problem). The proof of our pattern theorem for clusters is very different from either of the two methods for walks; our method does not work for walks, and the known walk methods do not work for clusters. Fortunately, the proof of Equation (??) does not depend much on the method of proof of the pattern theorem, and so we will be able to deduce (??) for weighted self-avoiding walks (see Section 3.5).

An interesting result of Bender, Gao, and Richmond (1992) is a graph-theoretical analogue of our pattern thereom. In particular, it shows that any fixed planar subgraph appears at least $\epsilon n$ times in almost all planar graphs with $n$ bonds. An important difference between the two papers is the lack of an underlying periodic lattice structure in Bender *et al.*

The pattern theorem has numerous consequences besides the ratio limit the-orem. Firstly, suppose that $L_1$ is a sublattice of $L$ (for example, the square lattice $\mathbf{Z}^2$ can be viewed as a sublattice of the triangular lattice; see Section 3.1). Clearly, the number of bond animals with $n$ sites on $L_1$ (up to translation) is bounded above by the number on $L$; hence the connective constant for bond animals on $L_1$ is less than or equal to that for $L$. This much is obvious, but the pattern theorem tells us that this inequality between growth constants is strict. See Section 3.4 for details. A similar result arises for some examples in which a set of clusters can be viewed as a subset of another set of clusters on the same lattice. For example, the pattern theorem provides a new proof of the result of Madras, Soteros, and Whittington (1988) that the growth constant for bond trees on $\mathbf{Z}^d$ is strictly smaller than that for bond animals. As a second example, Conway, Brak, and Guttmann (1993) present numerical estimates on several directed lattices that indicate that the growth constant for bond trees is strictly less than the growth constant for bond animals and strictly greater than the growth constant for site animals. In Sections 3.4 and 3.6 we prove these inequalities rigorously, as consequences of the pattern theorem. Further applications of the pattern theorem for weighted animals can be found in Sections 2.4 and 2.8 of Janse van Rensburg and Madras (1997). Some applications for self-avoiding walks are listed on page 231 of Madras and Slade (1993).

The remainder of the paper is organized as follows. The technical definitions and the statements of the main results are presented in Section 2. The assumptions for the results are phrased in terms of five Cluster Axioms that present five properties that are easy enough to check in most cases. The presentation of Section 2 gets somewhat abstract because we want to treat several different kinds of clusters, as well as general lattices. To help the reader, Section 3 presents many examples of the material of Section 2, arranged in parallel sections (e.g. Section 2.2 gives the abstract properties of clusters, while Section 3.2 gives examples of different families of clusters). In addition, Section 3.6 discusses modifications that need to be made for directed clusters. Section 4 presents the proofs of the two main theorems. Section 5 gives a summary and a discussion of some open problems.

## 2 Definitions and Results

This section will deal with definitions and statements of results for general undirected lattices. It is recommended that the reader refer to Section ?? for examples while reading the present section.

We shall view a $d$-dimensional lattice $L$ as a periodic embedding of an infinite graph in $\mathbf{R}^d$. The set of sites $S(L)$ is a countable subset of $\mathbf{R}^d$, and the set of bonds $B(L)$ is a countable set of unordered pairs of distinct sites (the sites are the "endpoints" of the bond). We write $\langle x, y \rangle$ to denote the bond whose endpoints are the sites $x$ and $y$. We assume two "local finiteness" conditions: each site is the endpoint of finitely many bonds, and that no bounded subset of $\mathbf{R}^d$ contains infinitely many sites. We also assume that $L$ is connected.

It will often be convenient to consider $L$ as a set of bonds and sites; that is, we identify $L$ with $S(L) \cup B(L)$. If $G$ is a subgraph of $L$ (or more generally a subset of $S(L) \cup B(L)$), then we write $S(G)$ and $B(G)$ respectively to denote the set of sites and the set of bonds of $G$.

**2.1 Translation:**

Let $G$ be a subgraph (or more generally a subset) of $L$, and let $u$ be a vector in $\mathbf{R}^d$. We write $G + u$ to denote the subgraph (or subset) which is the translation of $G$ by $u$:

$$S(G + u) = S(G) + u = \{x + u : x \in S(G)\}, \text{ and}$$

$$B(G + u) = B(G) + u = \{\langle x + u, y + u \rangle : \langle x, y \rangle \in B(G)\}.$$

Of course, this only makes sense if $S(G) + u \subset S(L)$ and $B(G) + u \subset B(L)$. This motivates us to define the set $S^*$ of all translations which leave $L$ invariant:

$$S^* = \{u \in \mathbf{R}^d : L + u = L\}.$$

Then $S^*$ is a group, and it can be shown that our "local finiteness" conditions ensure that $S^*$ is isomorphic to $\mathbf{Z}^k$ for some $k \leq d$ (Cassels (1959), Section III.4). We are interested in lattices with full-dimensional periodicity, so we shall assume that $S^*$ is isomorphic to $\mathbf{Z}^d$. (But see also example (j) in Section 3.1.) By applying a translation to $L$ if necessary, there is no loss of generality in assuming that the origin of $\mathbf{R}^d$, 0, is a site of $L$; indeed, we shall assume this throughout the paper. Now, since $0 \in S(L)$, it follows that $S^* \subset S(L)$. In many cases $S^* = S(L)$, but in other cases equality is not true. In general, the local finiteness conditions show that there is a finite set of sites, $a_1, a_2, \ldots, a_J$ (with $a_1 = 0$) such that $S^* + a_1, \ldots, S^* + a_J$ is a partition of $S(L)$. (That is, $\bigcup_{i=1}^{J}(S^* + a_i) = S(L)$, and $(S^* + a_i) \cap (S^* + a_j) = \emptyset$ whenever $i \neq j$).

For the remainder of this paper, the term "translation" will always mean "translation by an element of $S^*$".

Examples illustrating the preceding definitions may be found in Section 3.1.

**2.2 Clusters and weights:**

For each positive integer $n$, let $C_n$ be the set of all clusters of size $n$. In Section 3.2 we will give examples of what "cluster" and "size" could mean, but abstractly we only require that some Cluster Axioms be satisfied.

The first Cluster Axiom is simply:

> (*CA1*): $C_n$ is a collection of finite subgraphs of $L$ that is invariant under translation. (I.e., if $G \in C_n$ and $u \in S^*$, then $G + u \in C_n$.) The $C_n$'s are pairwise disjoint (i.e., $C_n \cap C_m = \emptyset$ whenever $n \neq m$).

Given (*CA1*), we can define $C_n^*$ to be the set of clusters in $C_n$ whose lexicographically smallest site is in $\{a_1, \ldots, a_J\}$. (We say that the point $(x_1, \ldots, x_d)$ is lexicographically smaller than $(y_1, \ldots, y_d)$ if there exists an $i \in \{1, \ldots, d\}$ such that $x_j = y_j$ for $j = 1, \ldots, i-1$ and $x_i < y_i$.) Thus, for each $G \in C_n$, there is a unique $H \in C_n^*$ and a unique $u \in S^*$ such that $G = H + u$. (If we think of translation as defining an equivalence relation among clusters, then we can view $C_n^*$ as the set of equivalence classes in $C_n$.) We also define

$$C_{<\infty} = \bigcup_{n=1}^{\infty} C_n$$

to be the set of all clusters.

We also define a weight function that assigns positive weight to each cluster

$$wt : C_{<\infty} \to (0, \infty)$$

that is invariant under translation:

$$wt(G) = wt(G + u) \quad \text{for every } u \in S^* \text{ and } G \in C_{<\infty}.$$

Our weight function should not be completely arbitrary. We shall need to know that changing a few sites and bonds in a cluster will only affect its weight within a bounded factor:

(*CA2*): For each $m$, there is a finite positive constant $\gamma_m$ with the property that

$$\frac{1}{\gamma_m} wt(G) \leq wt(G') \leq \gamma_m wt(G)$$

whenever $G$ and $G'$ differ by at most $m$ sites and bonds (i.e., whenever $|B(G) \Delta B(G')| + |S(G) \Delta S(G')| \leq m$, where $\Delta$ denotes symmetric difference).

Examples of weight functions satisfying (*CA2*) include the constant function ($wt(G) = 1$ for every $G$), as well as the collapse and percolation weights (recall (??) and (??)). A closely related weight function for animals is used in Madras et al. (1990): $wt(G) = z^{cyc(G)}$, where $z > 0$ and $cyc(G) = |B(G)| - |S(G)| + 1$ is the number of independent cycles in the animal $G$.

If $A$ is a subset of $C_{<\infty}$, then we write the weighted sum of members of $A$ as

$$\mathcal{G}(A) = \sum_{G \in A} wt(G). \tag{11}$$

We write $\mathcal{G}_n$ to denote the weighted sum of all clusters of size $n$ (up to translation),

$$\mathcal{G}_n = \mathcal{G}(C_n^*), \tag{12}$$

and we define
$$\lambda = \limsup_{n \to \infty}(\mathcal{G}_n)^{1/n}. \tag{13}$$

Our third Cluster Axiom is

(*CA3*): The limit $\lim_{n\to\infty}(\mathcal{G}_n)^{1/n}$ exists and is finite (and equals $\lambda$).

This axiom is known to hold for our main models of interest (see Section 3.2).

**2.3 Patterns:**

Let $P_1$ and $P_2$ be two finite disjoint subsets of $L$, with $P_1$ nonempty. Thus each $P_i$ can consist of bonds or edges or both (or neither, if $i = 2$); in particular, $P_i$ need not be a subgraph. Then the ordered pair $P = (P_1, P_2)$ is a *pattern*. If $G$ is a cluster, then we say that "$G$ contains $P$" if $G$ contains all of $P_1$ and none of $P_2$ (i.e., $P_1 \subset [B(G) \cup S(G)]$ and $P_2 \cap [B(G) \cup S(G)] = \emptyset$). If $x \in S^*$, then the translate of $P$ by $x$ is the new pattern

$$P + x = (P_1 + x, P_2 + x).$$

Thus the cluster $G$ contains $P + x$ if and only if $G - x$ contains $P$. We say that $P$ is a *proper pattern* if there are infinitely many values of $n$ such that $P$ is contained in some cluster of size $n$. This excludes patterns in which $P_2$ completely surrounds $P_1$, for example.

Our fourth Cluster Axiom says that any part of any cluster can be locally changed to create an occurrence of some translate of a given proper pattern.

> (*CA4*): For every proper pattern $P = (P_1, P_2)$, there exists a finite set $D$ of sites and bonds of $L$ (i.e., $D \subset S(L) \cup B(L)$) with the following property:
>
> For every cluster $G \in C_{<\infty}$ and every site $y \in S(G)$, there is another cluster $G'$ (possibly of different size) and a translation vector $t = t(y) \in S^*$ such that $y \in D + t$, $G'$ contains $P + t$, and $G' \setminus (D + t) = G \setminus (D + t)$.

(See Figure **??**.) We use the notation $G' = T(G, y)$. That is, $T$ is a function on clusters $G$ and their sites $y$ that creates a new cluster by altering sites and bonds inside a set $D + t$ around the specified site $y$ to create an occurrence of the pattern $P + t$, while leaving everything outside of $D + t$ the same. (This notation suppresses dependence on $P$.) In Section 3.3 we prove that animals and bond trees satisfy Axiom (*CA4*). The axiom must be modified slightly for directed clusters; see Section 3.6. We remark that self-avoiding walks do not satisfy Axiom (*CA4*).

In (*CA4*), the size of $D$ limits how different the sizes of $G$ and $G'$ can be. In particular, there exists a positive integer $\kappa$ (depending on $D$ and hence on $P$)

such that

$$T(G,y) \in \bigcup_{m=n-\kappa}^{n+\kappa} C_m \quad \text{whenever } G \in C_n \text{ and } y \in S(G). \tag{14}$$

Given Axiom $(CA4)$, Axiom $(CA2)$ now implies the existence of a constant $\gamma > 0$ (depending on $D$ and hence on $P$) such that

$$\frac{1}{\gamma} wt(G) \leq wt(T(G,y)) \leq \gamma \, wt(G) \quad \text{whenever } G \in C_{<\infty} \text{ and } y \in S(G). \tag{15}$$

### 2.4 The Pattern Theorem:

The first main theorem says that translates of a given proper pattern occur many times on most large clusters. More precisely, except for an exponentially small set of clusters of size $n$, the number of occurrences of translates of a given proper pattern is of the order $n$. Here, "exponentially small" is in terms of the weights of the clusters, not just the number of clusters.

**Theorem 2.1** *Assume that Cluster Axioms (CA1), (CA2), and (CA4) hold. Let $P$ be a proper pattern. Let $\mathcal{G}_n[\leq m, P]$ be the weighted sum of the set of clusters in $C_n^*$ which contain at most $m$ translates of $P$. Then there exists an $\epsilon > 0$ such that*

$$\limsup_{n \to \infty} (\mathcal{G}_n[\leq \epsilon n, P])^{1/n} < \lambda. \tag{16}$$

The proof is in Section 4. The important part of this result is that the inequality is strict. Recall that $\lambda$ was defined in Equation (??), and that if $(CA3)$ holds then $\lambda = \lim_{n \to \infty} (\mathcal{G}_n)^{1/n}$.

For some applications of this result, see Section 3.4.

### 2.5 Ratio Limit Theorem:

The ratio limit theorem requires some additional geometric information about the lattice, but as we shall see these conditions are easily verified for the examples mentioned in this paper.

Firstly, we need a special pair of patterns, which we shall call $U$ and $V$. Figure ?? shows one choice of the pair $U$ and $V$ for clusters in $\mathbf{Z}^2$. The important feature of these patterns is that any translate of $U$ in a cluster can be changed into a translate of $V$, with the size of the cluster increased by 1. (Also, the weight of the cluster changes by a precise multiplicative factor, which we call $\theta$.) Similarly, any $V$ can be changed easily into a $U$, with a decrease of 1 in cluster size. Now, the point of the Ratio Limit Theorem is to show that the sequence $\mathcal{G}_{n+1}/\mathcal{G}_n$ converges; once we know this, the value of the limit is obvious from $(CA3)$. So we want to show that $\mathcal{G}_{n+2}/\mathcal{G}_{n+1} \approx \mathcal{G}_{n+1}/\mathcal{G}_n$ when $n$ is large. By the Pattern Theorem, most large clusters contain lots of $U$'s and lots of $V$'s.

We could change any $U$ to a $V$, or any $V$ to $U$; thus typical clusters of size $n$ look very much like typical clusters of size $n + 1$ (or $n + 2$). The proof of the Ratio Limit Theorem works by considering all possible changes of one or two $U$'s into $V$'s in clusters of size $n$, and counting and comparing the results.

To give a formal description of the essential properties of $U$ and $V$, we shall introduce the following notation: For a pattern $P$ and a cluster $G$, define

$$\tau_P(G) := \{x \in S^* : G \text{ contains } P + x\}.$$

Then $|\tau_P(G)|$ is the number of translates of $P$ that occur in $G$. Also, if $\zeta$ is one of $U$ or $V$, then let $\hat{\zeta}$ be the other one (i.e., $\hat{\zeta} = V$ if $\zeta = U$, and $\hat{\zeta} = U$ if $\zeta = V$).

> ($CA5$) There exist proper patterns $U$ and $V$ and a constant $\theta \in (0, \infty)$, such that $U_1 \cup U_2 = V_1 \cup V_2$, and such that assertions ($i$) through ($iv$) hold whenever $\zeta \in \{U, V\}$ and $G$ is a cluster that contains $\zeta + x$ (i.e., $x \in \tau_\zeta(G)$):
>
> ($i$) $[G \setminus (\zeta_1 + x)] \cup (\hat{\zeta}_1 + x)$ is also a cluster, which we shall denote $\hat{G}_x$;
>
> ($ii$) $|\tau_\zeta(\hat{G}_x)| = |\tau_\zeta(G)| - 1$ and $|\tau_{\hat{\zeta}}(\hat{G}_x)| = |\tau_{\hat{\zeta}}(G)| + 1$;
>
> ($iii$) If $\zeta = U$ and $G \in C_n$, then $\hat{G}_x \in C_{n+1}$; and
>
> ($iv$) If $\zeta = U$, then $wt(\hat{G}_x) = \theta \, wt(G)$.

Informally, $\hat{G}_x$ is the result of changing one occurrence of $\zeta$ into $\hat{\zeta}$. Observe that if $\hat{G}_x = H$, then $\hat{H}_x = G$; hence if $\zeta = V$ and $G \in C_n$, then $\hat{G}_x \in C_{n-1}$ and $wt(\hat{G}_x) = \theta^{-1} \, wt(G)$.

**Theorem 2.2** *Assume that (CA1), (CA3), (CA5), and the conclusions of Theorem ?? all hold. (This will happen if we assume that all five Cluster Axioms hold.) Also assume that there exists a constant $\Upsilon$ such that $\mathcal{G}_{n+1} \geq \Upsilon \mathcal{G}_n$ for all sufficiently large $n$. Then*

$$\lim_{n \to \infty} \frac{\mathcal{G}_{n+1}}{\mathcal{G}_n} = \lambda.$$

The proof appears in Section 4. The verification of the assumptions of Theorem ?? for the main models, including the proof of Theorem ??, is in Section 3.5.

## 3 Examples

In this section we shall illustrate the definitions and results of Section ?? by various examples.

**3.1 Lattices:**

We list several standard lattices, as well as a few non-standard ones. Figure ?? illustrates some of them.

(a) <u>Hypercubic ($Hyp_d$)</u>: The $d$-dimensional hypercubic lattice is the lattice whose sites are the points of $\mathbf{Z}^d$ and whose bonds join nearest-neighbour pairs. As is customary, we shall often denote this lattice by $\mathbf{Z}^d$. This ambiguity of $\mathbf{Z}^d$ representing a lattice as well as a set of sites should not lead to confusion in the rest of the paper, but for precision in this section we shall use $Hyp_d$ to denote the $d$-dimensional hypercubic lattice. Thus, we have

$$S(Hyp_d) = \mathbf{Z}^d \quad \text{and} \quad B(Hyp_d) = \{\langle x, y \rangle : x, y \in \mathbf{Z}^d, ||x - y||_1 = 1\}$$

where $||(u_1, \ldots, u_d)||_1 = |u_1| + \cdots + |u_d|$. For this lattice we have $S^* = \mathbf{Z}^d$, and therefore $J$ equals 1 and $a_1$ is the origin.

(b) <u>Triangular lattice ($Tri$)</u>: This two-dimensional lattice can be represented by

$$S(Tri) = \mathbf{Z}^2 \quad \text{and} \quad B(Tri) = B(Hyp_d) \cup \{\langle x, x + (1, 1) \rangle : x \in \mathbf{Z}^2\}.$$

Then $S^* = \mathbf{Z}^2$, $J = 1$ and $a_1 = 0$.

(c) <u>Hexagonal lattice ($Hex$)</u>: We can represent this lattice as a sublattice of $Hyp_2$ as follows:

$$S(Hex) = \mathbf{Z}^2 \quad \text{and}$$

$$B(Hex) = B(Hyp_2) \setminus \{\langle (a, b), (a, b + 1) \rangle : a, b \in \mathbf{Z},\ a + b \text{ is odd}\}.$$

For this lattice we have $S^* = \{(a, b) \in \mathbf{Z}^2 : a + b \text{ is even}\}$ and $J = 2$, with $a_1 = (0, 0)$ and $a_2 = (1, 0)$.

(d) <u>Kagome lattice ($Kag$)</u>: To represent this two-dimensional lattice, we shall write $2\mathbf{Z}^2$ to denote the points of $\mathbf{Z}^2$ having both coordinates even. Let $a_1 = (0, 0)$, $a_2 = (1, 0)$, and $a_3 = (0, 1)$. Then we have

$$S(Kag) = \cup_{i=1}^{3}(a_i + 2\mathbf{Z}^2) \quad \text{and}$$

$$B(Kag) = \{\langle a_1, a_2 \rangle, \langle a_1, a_3 \rangle, \langle a_2, a_3 \rangle, \langle a_2, (2, 0) \rangle, \langle a_3, (0, 2) \rangle, \langle a_3, (-1, 2) \rangle\} + 2\mathbf{Z}^2.$$

Then $S^* = 2\mathbf{Z}^2$, $J = 3$, and the $a_i$'s are as given above.

(e) <u>Rectangular $(r_1, r_2)$ lattice ($Rect_{r_1, r_2}$)</u>: In this family of two-dimensional lattices, $r_1$ and $r_2$ could be any positive integers. We can describe $Rect_{r_1, r_2}$ as the intersection of $Hyp_2$ with all lines of the form $x_i = kr_i$ ($k \in \mathbf{Z}$; $i = 1, 2$). More formally,

$$S(Rect_{r_1, r_2}) = \{(kr_1, b) : k, b \in \mathbf{Z}\} \cup \{(a, kr_2) : a, k \in \mathbf{Z}\} \quad \text{and}$$
$$B(Rect_{r_1, r_2}) =$$
$$\{\langle (kr_1, b), (kr_1, b + 1) \rangle : k, b \in \mathbf{Z}\} \cup \{\langle (a, kr_2), (a + 1, kr_2) \rangle : a, k \in \mathbf{Z}\}.$$

Then $S^* = \{(k_1 r_1, k_2 r_2) : k_1, k_2 \in \mathbf{Z}\}$ and $J = r_1 + r_2 - 1$.

(f) <u>$d$-dimensional spread-out lattice of range $M$ ($\mathbf{Z}^d_{(M)}$)</u>: Let $M$ be a positive real number and let $||\cdot||$ be a norm on $\mathbf{R}^d$. Then $S(\mathbf{Z}^d_{(M)}) = \mathbf{Z}^d$ and
$$B(\mathbf{Z}^d_{(M)}) = \{\langle x, y\rangle : x, y \in \mathbf{Z}^d, 0 < ||x - y|| \leq M\}.$$
Here, $S^* = \mathbf{Z}^d$. These lattices have been used to approximate "mean field" behaviour, often with the sup norm $||x||_\infty = \max\{|x_1|, \ldots, |x_d|\}$ (e.g. Hara and Slade (1992)).

(g) <u>Dead-end lattice ($DE$)</u>: To describe this unusual two-dimensional lattice, let $a_1 = (0,0)$ and $a_2 = (1/2, 1/2)$. Let
$$S(DE) = \mathbf{Z}^2 \cup (a_2 + \mathbf{Z}^2) \quad \text{and}$$
$$B(DE) = B(Hyp_2) \cup (\langle a_1, a_2\rangle + \mathbf{Z}^2).$$
Then $S^* = \mathbf{Z}^2$ and $J = 2$. This lattice and similar ones serve as a class of counterexamples, but could also be used to model clusters of $Hyp_2$ in which sites can be of two types.

(h) <u>Body-centered cubic lattice ($BCC$)</u>: This classical lattice in $\mathbf{R}^3$ has
$$S(BCC) = \{(x_1, x_2, x_3) \in \mathbf{Z}^3 : x_1 + x_2 + x_3 \text{ is a multiple of 3}\},$$
$$B(BCC) = \{\langle x, y\rangle : |x_1 - y_1| = |x_2 - y_2| = |x_3 - y_3| = 1\},$$
and $S^* = S(BCC)$.

(i) <u>Face-centered cubic lattice ($FCC$)</u>: This classical lattice in $\mathbf{R}^3$ has
$$S(FCC) = \{(x_1, x_2, x_3) \in \mathbf{Z}^3 : x_1 + x_2 + x_3 \text{ is a multiple of 2}\},$$
$$B(FCC) = \{\langle x, y\rangle : x, y \in \mathbf{Z}^3, |x_1 - y_1|^2 + |x_2 - y_2|^2 + |x_3 - y_3|^2 = 2\},$$
and $S^* = S(FCC)$.

(j) <u>Slabs</u>: A $k$-dimensional slab of the $d$-dimensional lattice $L$ is the part of $L$ that lies between $d - k$ given pairs of parallel hyperplanes. (See Section 6.4 of Grimmett (1989) or Section 8.2 of Madras and Slade (1993) for some problems related to slabs.) For example, if $M$ and $M'$ are positive integers, then
$$Hyp_d \cap \{(z_1, \ldots, z_d) \in \mathbf{R}^d : 0 \leq z_{d-1} \leq M, 0 \leq z_d \leq M'\}$$
is a $(d-2)$-dimensional slab of $Hyp_d$. We can view it as a $(d-2)$-dimensional lattice via the following mapping from $\mathbf{R}^d$ to $\mathbf{R}^{d-2}$, which is one-to-one on the sites of the slab in $\mathbf{Z}^d$:
$$(z_1, \ldots, z_d) \mapsto \left(z_1 + \frac{z_{d-1}}{M+1} + \frac{z_d}{(M+1)(M'+1)}, z_2, \ldots, z_{d-2}\right).$$
The image of this slab in $\mathbf{R}^{d-2}$ is a lattice with $S^* = \mathbf{Z}^{d-2}$.

Other unusual lattices may be found in Conway, Brak, and Guttmann (1993).

**3.2 Clusters and weights:**

The following are examples of sets that could be considered for $C_n$, the set of clusters of size $n$ on a lattice $L$. (Note: We are not claiming that all of them satisfy all of the cluster axioms. This will be discussed later.)

(a) The set of bond animals of $L$ which contain exactly $n$ sites. Recall that a bond animal is a finite connected subgraph of the infinite graph $L$.

(a$'$) The set of bond animals of $L$ which contain exactly $n$ bonds.

(b) The set of bond trees of $L$ which contain exactly $n$ sites. Recall that a bond tree is a bond animal with no cycles; thus every tree with $n$ sites has $n-1$ bonds.

(c) The set of site animals of $L$ which contain exactly $n$ sites. A site animal is a finite connected subgraph $G$ of $L$ whose bonds are determined by its sites in the sense that
$$B(G) = \{\langle x,y \rangle \in B(L) : x, y \in S(G)\}.$$

(d) The set of directed bond animals of $L$ which contain exactly $n$ sites. To describe these objects, fix a nonzero vector $\vec{v} \in \mathbf{R}^d$ such that $\vec{v} \cdot x \neq \vec{v} \cdot y$ for every bond $\langle x, y \rangle$ of $L$ (i.e., $\vec{v}$ is not orthogonal to any bond). Suppose that $G$ is a subgraph of $L$, and that $y$ and $z$ are two sites of $G$. We say that there is a $\vec{v}$-directed path from $y$ to $z$ in $G$ if there is a finite sequence of sites $(x^{(i)} : i = 0, \ldots, k)$ in $G$ such that $x^{(0)} = y$ and $x^{(k)} = z$, and $\langle x^{(i)}, x^{(i+1)} \rangle \in B(G)$ and $\vec{v} \cdot x^{(i)} < \vec{v} \cdot x^{(i+1)}$ for each $i = 0, \ldots, k-1$. (Here $\vec{v} \cdot x$ is the usual Euclidean inner product in $\mathbf{R}^d$.) For example, if $L = Hyp_3$ and $\vec{v} = (1,1,1)$, then the $\vec{v}$-directed paths are those that only take steps in the positive coordinate directions. A *$\vec{v}$-directed bond animal* is a bond animal that contains a site $r$ with the property that there are $\vec{v}$-directed paths from $r$ to every other site of the animal. The site $r$ is called the *root* of the animal. We will often omit the prefix $\vec{v}$ in our terminology.

(e) The set of directed bond trees with $n$ sites. A $\vec{v}$-directed bond tree is a bond tree (in the undirected sense of (b) above) that is also a $\vec{v}$-directed bond animal. Equivalently, a $\vec{v}$-directed bond tree is a $\vec{v}$-directed bond animal in which every site $x$ except the root has exactly one "incoming" bond (i.e., a bond $\langle w, x \rangle$ such that $\vec{v} \cdot w < \vec{v} \cdot x$) in the animal.

(f) The set of self-avoiding polygons in $Hyp_d$ with $2n$ bonds. A self-avoiding polygon is a bond animal $G$ in which every site of $G$ is the endpoint of exactly two bonds of $G$.

(g) The set of self-avoiding walks in $Hyp_d$ with $2n$ bonds. A self-avoiding walk is a bond tree $G$ in which no site of $G$ is the endpoint of more than two bonds of $G$. (We use $2n$ instead of $n$ here so that Axiom $(CA5)$ will hold.)

(g$'$) The set of self-avoiding walks in $Hyp_d$ with $2n+1$ bonds.

Other examples are evident: Site animals containing $n$ bonds; directed site animals containing $n$ sites; etc.

Axiom ($CA3$) holds for all of the above examples with weights as in (??), thanks to concatenation and subadditivity arguments; see Klarner (1967), Klein (1981), Soteros and Whittington (1990), Madras *et al.* (1990), Madras and Slade (1993), and Janse van Rensburg and Madras (1997).

**3.3 Patterns:**

We begin with some remarks about proper patterns. Let $L$ be a lattice. If $P_2$ is the empty set, then $(P_1, \emptyset)$ is a proper pattern for bond animals or site animals for every set $P_1$; but this need not be true for all kinds of clusters. Obviously, if $P_1$ contains a cycle, then $(P_1, P_2)$ cannot be a proper pattern for trees. If $P_1$ is a subgraph of $L$ and $P_2$ includes all bonds of $L$ that have exactly one endpoint in $P_1$, then $(P_1, P_2)$ cannot be a proper pattern for any class of connected clusters, since no cluster containing $P$ could contain any site outside $P_1$.

Next we shall prove that the fourth Cluster Axiom holds in a wide range of cases.

**Proposition 3.1** *Let $L$ be any of the lattices described in Section 3.1. Then Axiom (CA4) holds for bond animals, site animals, and bond trees.*

**Proof**: We will first prove the result for bond (or site) animals on the hypercubic lattice $Hyp_d$. We then outline the extension to other lattices of Section 3.1. Finally, we describe the proof for bond trees.

Let $P = (P_1, P_2)$ be a proper pattern. Choose an integer $M$ such that $P_1 \cup P_2 \subset \{x \in \mathbf{R}^d : ||x||_\infty \leq M\}$ (where $||x||_\infty = \max\{|x_1|, \ldots, |x_d|\}$), and we use the natural convention that a bond $\langle u, v \rangle$ is contained in a set if and only if the set contains the line segment joining $u$ and $v$). Let $D$ be the set of all sites and bonds of $L$ in the cube $\{x \in \mathbf{R}^d : ||x||_\infty \leq M + 1\}$, and let $\partial D$ be the set of all sites and bonds of $L$ in $\{x \in \mathbf{R}^d : ||x||_\infty = M + 1\}$, the boundary of $D$.

Let $H$ be a bond (or site) animal that contains $P$ and has at least one site outside $D$ (we can do this because $P$ is proper), and let $\tilde{H} = (H \cap D) \cup \partial D$. (See Figure ??). Then $\tilde{H}$ is an animal that contains $P$ (in particular, it is connected because $\partial D$ is connected and every path in the lattice from $D$ to $D^c$ contains a site of $\partial D$). As in Axiom ($CA4$), suppose we are given an animal $G$ and a site $y \in S(G)$. Define $t(y) = y$ and

$$G' = T(G, y) = \big(G \setminus (D + t)\big) \cup (\tilde{H} + t).$$

The picture is that $G'$ contains the "surface" $\partial D$ that has been translated to surround $y$, agrees with $G$ outside this surface, and looks like $H$ inside this surface. (See Figure ??(c,d).) It is not hard to see that this produces an animal with the desired properties. Thus Axiom ($CA4$) holds for bond and site animals.

For other lattices, we can also choose a set $D$ of the form $L \cap \{x \in \mathbf{R}^d : ||x||_\infty \leq M_1\}$ and a set $\partial D$ of the form $L \cap \{x \in \mathbf{R}^d : M_0 \leq ||x||_\infty \leq M_1\}$ for suitably chosen $M_1$ and $M_0$. The key properties that guide this choice are: $(i)$ $P_1 \cup P_2 \subset D \setminus \partial D$; $(ii)$ $\partial D$ is connected; $(iii)$ every path in the lattice from $D$ to $D^c$ contains a site of $\partial D$; and $(iv)$ $D \cap (S^* + a_i) \neq \emptyset$ for every $i = 1, \ldots, J$. For example, for the spread-out lattice $\mathbf{Z}_{(M)}^d$ with norm $||\cdot||_\infty$, choose $M_0 > \max\{||x||_\infty : x \in P_1 \cup P_2\}$ and $M_1 = M_0 + M$. Next, for the translation vectors $t$: Given a site $y$, choose $t = t(y) \in S^*$ so that $y \in D + t(y)$ (this can be done by property $(iv)$). The proof now proceeds as for the hypercubic case.

Finally, consider the case of bond trees. Let $D$ and $\partial D$ be as above. Let $H$ be a bond tree that contains $P$ and has at least one site outside $D$, and let $\tilde{H} = (H \cap D) \cup \partial D$. Then $\tilde{H}$ is a bond animal that contains $P$, but in general it is not a tree. For Axiom $(CA4)$, suppose we are given a tree $G$ and a site $y \in S(G)$. Let $t = t(y)$ be as above, and let

$$G_A = \bigl(G \setminus (D + t)\bigr) \cup (\tilde{H} + t);$$

then $G_A$ is a bond animal that contains $P + t$. Next, let $G_B$ be the subgraph of $G_A$ obtained by deleting all bonds in $D$ that have at least one endpoint in $(\partial D) + t$. Then $G_B$ contains $P + t$, but it is disconnected. However, $G_B$ contains no cycles (since $G_B$ is the disjoint union of a subgraph of the tree $G$, a subgraph of $(H \cap D) + t$, and possibly some isolated sites of $(\partial D) + t$). The existence of the graph $G'$ in $(CA4)$ is now guaranteed by the following routine exercise of graph theory: *Let $G_B$ be subgraph of a connected graph $G_A$. If $G_B$ contains no cycles, then $G_A$ contains a tree which contains $G_B$.* □

### 3.4 The Pattern Theorem:

The preceding parts of this section have shown that the assumptions of the Pattern Theorem ?? hold for bond animals, site animals, and bond trees on any of the lattices of Section 3.1, with any weights satisfying $(CA2)$. Section 3.6 discusses the situation for directed clusters.

The Pattern Theorem implies strict inequality between $\lambda$'s for different families of clusters. In this subsection, we shall use the notation $\lambda_{BA}[L]$ to denote the value of $\lambda$ for bond animals on the lattice $L$ (we suppress notational dependence on the choice of weights). We replace $BA$ by $BT$ for bond trees, and by $SA$ for site animals. Some of the following results apply only to cluster weights that are identically 1, i.e. to the case $\mathcal{G}_n = |C_n^*|$. In this case $\lambda$ is defined by Equation (??) and is called the growth constant. We shall write $\mu$ instead of $\lambda$ for the growth constant; e.g. $\mu_{SA}[L]$ denotes the growth constant for site animals on $L$.

**Corollary 3.2** *Let $L_1$ be a sublattice of $L_2$ (i.e., $S(L_1) \subset S(L_2)$, $B(L_1) \subset B(L_2)$, and $L_1 \neq L_2$), both satisfying the properties of Section 2.1. Then $\mu_{BA}[L_1] < \mu_{BA}[L_2]$, $\mu_{BT}[L_1] < \mu_{BT}[L_2]$, and $\mu_{SA}[L_1] < \mu_{SA}[L_2]$. (Here we could be measuring cluster size either by number of sites or by number of bonds.)*

**Proof:** First we consider bond animals and bond trees. Fix a bond $b \in B(L_2) \setminus B(L_1)$. Then $P = (\{b\}, \emptyset)$ is proper pattern for clusters on $L_2$, but never occurs in clusters of $L_1$. The corollary is thus an immediate consequence of the Pattern Theorem.

Now consider site animals. Notice that a site animal on $L_1$ need not be a site animal on $L_2$ (e.g. a unit square of four sites and four bonds is a site animal on the square lattice, but not on the triangular lattice (Section 3.1(b)), since it is missing the diagonal bond). However, identifying a site animal with its set of sites, it is clear that every site animal in $L_1$ corresponds to a site animal in $L_2$ with the same set of sites (and this correspondence is one-to-one, but not onto). The proof of the preceding paragraph works for site animals if $S(L_1) \neq S(L_2)$ (using a site for $P_1$ instead of a bond), so assume that $S(L_1) = S(L_2)$. Again, fix a bond $b \in B(L_2) \setminus B(L_1)$. Since $L_2$ is infinite and connected, we can choose a sequence $(x^{(0)}, x^{(1)}, \ldots)$ of distinct sites of $L_2$ such that $b = \langle x^{(0)}, x^{(1)} \rangle$, $\langle x^{(i)}, x^{(i+1)} \rangle \in B(L)$ for every $i \geq 0$ and $\langle x^{(0)}, x^{(j)} \rangle \notin B(L)$ for every $j \geq 2$. Let $N$ be the set of those sites of $L_2$ which are neighbours of $x^{(0)}$, except for $x^{(1)}$:

$$N = \{z \in S(L_2) : \langle x^{(0)}, z \rangle \in B(L), z \neq x^{(1)}\}.$$

Let $P_1 = \{x^{(0)}, x^{(1)}\}$ and $P_2 = N$. Then $P = (P_1, P_2)$ is a proper pattern for site animals in $L_2$, but cannot occur on any site animal in $L_1$ (since any large subgraph of $L_1$ containing $P$ cannot be connected). Therefore the inequality $\mu_{SA}[L_1] < \mu_{SA}[L_2]$ follows from the Pattern Theorem. □

We remark that the preceding proof does not apply to general weight functions, since the same cluster could have different weights on different lattices. For example, removing some bonds from a lattice can change the number of contact bonds. Similar things happen in part (ii) of the next result.

**Corollary 3.3** *Let $L$ be any of the lattices of Section 3.1. Assume that we measure the size of a cluster by the number of sites.*
*(i) For any weights satisfying (CA2), $\lambda_{BT}[L] < \lambda_{BA}[L]$.*
*(ii) For weights identically 1, $\mu_{SA}[L] < \mu_{BT}[L]$.*

**Proof:** (i) Let $P_1$ be a cycle of $L$ and let $P_2 = \emptyset$. Then $P = (P_1, \emptyset)$ is a proper pattern for bond animals on $L$. Since bond trees contain no translates of $P$, the Pattern Theorem says that they must be exponentially rare in the set of bond animals.

(ii) We shall use the notation $C^*_{SA,n}$ and $C^*_{BT,n}$ to distinguish the set of site animal clusters from the set of bond tree clusters. Define the map $\Phi : C^*_{BT,n} \to C^*_{SA,n}$ so that $\Phi(G)$ is the unique $G'$ in $C^*_{SA,n}$ such that $S(G') = S(G)$. That is, $\Phi$ fills in the "missing bonds" of the tree $G$. The map $\Phi$ is clearly onto.

As in part (i), let $P_1$ be a cycle of $L$ and let $P = (P_1, \emptyset)$. Then $P$ is a proper pattern for site animals. Next, let $\alpha$ and $\beta$ be two different bonds of $P_1$. For $i = \alpha$ or $\beta$, let $P_1^{(i)} = P_1 \setminus \{i\}$, $P_2^{(i)} = \{i\}$, and $P^{(i)} = (P_1^{(i)}, P_2^{(i)})$. Then $P^{(i)}$ is a proper pattern for bond trees, but cannot occur in site animals.

Observe that there exists a $K > 0$, depending on $P$, such that every cluster containing $m$ translates of $P$ must contain at least $m/K$ disjoint translates of $P$ (i.e., the corresponding translates of $P_1 \cup P_2$ are disjoint). For $\epsilon > 0$, let $C_{SA,n}^{disj}(\epsilon)$ be the set of all clusters in $C_{SA,n}^*$ that contain at least $\epsilon n$ disjoint translates of $P$. The preceding observation and the Pattern Theorem tells us that there is an $\epsilon > 0$ such that $|C_{SA,n}^{disj}(\epsilon)| > |C_{SA,n}^*|/2$ for all sufficiently large $n$.

Let $G'$ be an arbitrary site animal in $C_{SA,n}^{disj}(\epsilon)$. Then there are at least $2^{\epsilon n}$ trees $G$ in $C_{BT,n}^*$ such that $\Phi(G) = G'$. (This is because each translate of $P$ in $G'$ could have arisen from a translate of either $P^{(\alpha)}$ or $P^{(\beta)}$ in $G$). Hence

$$|C_{BT,n}^*| \geq 2^{\epsilon n} |C_{SA,n}^{disj}(\epsilon)| > 2^{\epsilon n - 1} |C_{SA,n}^*|$$

for all sufficiently large $n$. Taking $n^{th}$ roots and letting $n \to \infty$ shows that $\mu_{BT}[L] \geq 2^{\epsilon} \mu_{SA}[L]$, and the result follows. □

### 3.5 Ratio Limit Theorem:

First we shall show that Axiom $(CA5)$ holds in our main models of interest. The proof of Theorem 7.3.2 in Madras and Slade (1993) shows that it holds for self-avoiding walks and self-avoiding polygons on $\mathbf{Z}^d$. For animals and trees, we have the following result.

**Proposition 3.4** *Let $L$ be any of the lattices of Section 3.1. Consider weights of the form (??) for (bond or site) animals or bond trees. Then there exists a pair of patterns $U$ and $V$ that satisfy Axiom (CA5).*

**Proof:** Since $L$ is an infinite connected graph, there exists an infinite sequence of distinct sites $(x^{(0)}, x^{(1)}, \ldots)$ such that $\langle x^{(i)}, x^{(i+1)} \rangle \in B(L)$ for every $i \geq 0$ and $\langle x^{(0)}, x^{(j)} \rangle \notin B(L)$ for every $j \geq 2$. Let $N$ be the set of those sites of $L$ which are neighbours of $x^{(0)}$, except for $x^{(1)}$:

$$N = \{z \in S(L) : \langle x^{(0)}, z \rangle \in B(L), z \neq x^{(1)}\}.$$

Let $P_1 = \{x^{(1)}\}$ and and $P_2 = N \cup \{x^{(0)}, \langle x^{(0)}, x^{(1)} \rangle\}$. Then $P = (P_1, P_2)$ is a proper pattern (since for any $n \geq 1$, the sites $x^{(1)}, \ldots, x^{(n)}$ and the corresponding bonds form a cluster).

Given this $P$, choose $D$, $\partial D$, $H$, and $\tilde{H}$ as in the proof of Proposition ??. For the case of bond or site animals, define the first pattern $U$ by $U_1 = \tilde{H}$ and $U_2 = D \setminus \tilde{H}$. Define the second pattern $V$ by $V_1 = \tilde{H} \cup \{x^{(0)}, \langle x^{(0)}, x^{(1)} \rangle\}$ and $V_2 = D \setminus V_1$. (See Figure ??.) Observe that $U$ and $V$ are proper patterns for bond or site animals. (The $V$ for bond trees is an appropriate spanning tree of the $V$ for bond animals, as in the last paragraph of the proof of Proposition ??; then $U$ is obtained by deleting $\{x^{(0)}, \langle x^{(0)}, x^{(1)} \rangle\}$ from $V$.) For weights of the form (??), we can then take $\theta = z_s^{|N|-1}$ in $(CA5)(iv)$. To see that $(CA5)(ii)$ holds, notice that if two translates of $U$ or $V$ in a cluster $G$ overlap, then they must only overlap in the translates of $\partial D$; therefore changing a $U$ to a $V$ (say) cannot affect any other occurrence of $U$ or $V$ in the cluster. It is now routine to check that Axiom ($CA5$) holds for bond (and site) animals as well as bond trees. □

Next we shall show that the final assumption of Theorem ?? holds in a wide class of models.

**Proposition 3.5** *Assume that the lattice $L$ has the following property: Every site $x \in S(L)$ is the endpoint of a bond $\langle x, y \rangle \in B(L)$ such that $y$ is lexicographically greater than $x$. Then for all examples of <u>undirected</u> clusters from Section 3.2, with weights satisfying Axiom (CA2), there exists a constant $\Upsilon$ (depending on the model) such that $\mathcal{G}_{n+1} \geq \Upsilon \mathcal{G}_n$ for all sufficiently large $n$.*

**Remark**: The assumption of Proposition ?? holds for every lattice of Section 3.1, except for the Dead-End lattice ($DE$). However, the result of the proposition does hold for $DE$ because the assumption is true for the lattice obtained by reflecting $DE$ through the origin.

**Proof of Proposition ??**: First we consider animals and trees. Fix $n$. For an arbitrary cluster $G \in C_n^*$, let $x_G$ be the lexicographically largest site of $G$. By our assumption, there exists a bond $\langle x_G, y_G \rangle$ in $B(L)$ such that $y_G$ is lexicographically larger than $x_G$. For the cases of bond animals or bond trees, let $G^+$ be the cluster $G \cup \{y_G, \langle x_G, y_G \rangle\}$; for site animals (counted by sites), let $G^+$ be the cluster defined by $S(G^+) = S(G) \cup \{y_G\}$. Then $G^+ \in C_{n+1}^*$. In fact, the map $G \mapsto G^+$ is one-to-one (because the lexicographically largest site of $G^+$ must be $y_G$). Also, $wt(G^+) \geq wt(G)/\gamma_2$ by ($CA2$). The result follows with $\Upsilon = 1/\gamma_2$.

The argument for self-avoiding walks and self-avoiding polygons proceeds as on page 230 of Madras and Slade (1993). □

**Remark:** The above proof is not quite complete for the case of site animals counted by number of bonds, since $G^+$ as given may contain more than $n+1$ bonds. However, it is possible to prove Proposition ?? for this class of clusters for any of the lattices of Section 3.1 by considering each lattice separately.

As a consequence of the preceding results, we have the following.

**Corollary 3.6** *The Ratio Limit Theorem ?? holds for any of the undirected clusters of Section 3.2 on any of the lattices of Section 3.1, with weights of the form (??). In particular, Theorem ?? holds.*

In Section 3.6 we shall show that Theorem ?? also holds for directed clusters.

### 3.6 Directed clusters:

The definitions of directed paths, animals and trees were given in Section 3.2(d,e). We assume the following directed connectivity property of the lattice $L$: Whenever $i$ and $j$ are in $\{1, \ldots, J\}$, there exists a $\vec{v}$-directed path from $a_i$ to some site of $S^* + a_j$ (and hence, by translation invariance, there exists a $\vec{v}$-directed path from some site of $S^* + a_i$ to $a_j$). Note that this excludes unusual examples such as the Dead-End lattice.

A pattern $P = (P_1, P_2)$ is a *proper pattern* if for every finite subset $F$ of $L$, there exists a cluster $G$ containing $P$ whose root is outside $F$. As an example, consider $(1,1)$-directed bond animals in the square lattice $\text{Hyp}_2$. If $P_1$ contains the origin and $P_2$ contains the two sites $(-1, 0)$ and $(0, -1)$, then $(P_1, P_2)$ is not a proper pattern for these clusters because any such cluster that contains $(P_1, P_2)$ must have $(0, 0)$ as its root.

We shall use the following Directed Cluster Axiom $(DCA4)$ instead of $(CA4)$.

> $(DCA4)$: For every proper pattern $P = (P_1, P_2)$, there exist finite sets $D_1, \ldots, D_J$ of sites and bonds of $L$ (i.e., $D_i \subset S(L) \cup B(L)$) with the following property:
>
> For every cluster $G \in C_{<\infty}$ and every site $y \in S(G)$, there is another cluster $G'$ (possibly of different size) and a translation vector $t = t(y) \in S^*$ such that $y \in D_i + t$, $G'$ contains $P + t$, and $G' \setminus (D_i + t) = G \setminus (D_i + t)$ (where $i$ is the subscript such that $y \in S^* + a_i$).

Again, we use the notation $G' = T(G, y)$.

We have the following analogue of Proposition ??.

**Proposition 3.7** *Let $L$ be any of the lattices described in Section 3.1 (except the Dead-End lattice). Then Axiom (DCA4) holds for directed bond animals, directed site animals, and directed bond trees.*

**Proof**: The proof is similar to the proof of Proposition ??. For the case of $\text{Hyp}_d$, we take $D_1 = D = L \cap \{x \in \mathbf{R}^d : \|x\|_\infty \leq M + 1\}$, as defined in the proof of Proposition ??. We define $t(y) = y - z_D$, where $z_D$ is the "corner" site of $D$ such that $\vec{v} \cdot z_D < \vec{v} \cdot x$ for every other site of $D$. For directed animals on other lattices, we choose the sets $D$ and $\partial D$ satisfying the properties $(i)$ $P_1 \cup P_2 \subset D \setminus \partial D$; $(ii')$ there exist sites $z_D$ and $z^D$ in $S(D)$ with the property that for every $x \in S(\partial D)$, there exists a $\vec{v}$-directed path in $\partial D$ from $z_D$ to $z^D$ that contains $x$; and $(iii')$ every path from a site of $D$ to a site of $D^c$ (or vice versa) contains a site of $\partial D$. See Figure ??.

Let $H$ be a directed animal containing $P$ whose root is outside $D$. Let $\tilde{H} = (H \cap D) \cup \partial D$; this is a directed animal that contains $P$. For each $i = 1, \ldots, J$, we proceed as follows. By the directed connectivity assumption on $L$, there exists a site $z_{(i)} \in S^* + a_i$ such that there is a $\vec{v}$-directed path from $z_{(i)}$ to $z_D$. Let $\tilde{H}_i$ be the union of $\tilde{H}$ and this path; also let $D_i$ be the union of $D$ and this path. Observe that $\tilde{H}_i$ and $D_i$ are directed animals rooted at $z_{(i)}$. Now, if $y \in S(G) \cap (S^* + a_i)$, then we let $t = t(y) = y - z_{(i)}$ and

$$G' \;=\; T(G, y) \;=\; (G \setminus (D_i + t)) \cup (\tilde{H}_i + t)\ .$$

Observe that $D_i + t$ and $\tilde{H}_i + t$ are rooted at $y$. With these constructions, the proof for directed animals is essentially the same as for the undirected case of Proposition ??.

The case of directed bond trees also is similar to the proof for bond trees in Proposition ??. In particular, we use the following: Let $G_A$ be a $\vec{v}$-directed bond animal with root $r$, and let $G_B$ be a subgraph of $G_A$ (in the undirected sense). If every site of $G_B$ has at most one incoming bond (as defined in Section 3.2(e)) in $G_B$, then $G_A$ contains a $\vec{v}$-directed bond tree that contains $G_B$ and has root $r$. (Proof: Let $\tau = G_B \cup \{r\}$. If $\tau$ is a directed bond tree, then we are done. If not, then there is site $w \in S(\tau) \setminus \{r\}$ which has no incoming bond in $\tau$. But $w$ has an incoming bond in $G_A$; so add this bond (and its other endpoint, if necessary) to $\tau$. Repeat this procedure with the new $\tau$. Continue until every site of $\tau$ except $r$ has an incoming bond. The final $\tau$ will be what we want.) $\square$

Theorem ?? holds for directed bond animals, site animals, and bond trees if we replace $(CA4)$ by $(DCA4)$; indeed, the proof of the theorem is the same (see Section ??). The directed analogues of the ensuing Corollaries ?? and ?? of Section 3.4 also hold. The proofs of these corollaries carry over, with the following modification to the proof of Corollary ??$(ii)$: There exist sites $y$ and $z$ of $L$ and $\vec{v}$-directed paths $\pi_1$ and $\pi_2$ from $y$ to $z$ such that $\pi_1 \cap \pi_2 = \{y, z\}$ (here each $\pi_i$ is a set of sites and bonds). Let $\alpha$ (respectively, $\beta$) be the bond of $\pi_1$ (respectively, $\pi_2$) that has $z$ for an endpoint. Let $P_1 = \pi_1 \cup \pi_2$. With these definitions of $P_1$, $\alpha$, and $\beta$, the rest of the proof is unchanged.

The directed analogue of Proposition ??, showing that Axiom $(CA5)$ holds for directed animals and bond trees, can be proven with the following modification in the definition of the $x^{(i)}$'s. Let $x^{(0)}$ be an arbitrary site of $L$. For $i = 0, 1, \ldots$, inductively define $x^{(i+1)}$ to be the neighbour $w$ of $x^{(i)}$ that minimizes $\vec{v} \cdot w$. By the directed connectivity assumption on $L$, we know that $\vec{v} \cdot x^{(i+1)} < \vec{v} \cdot x^{(i)}$ for every $i$. Also observe that $\langle x^{(0)}, x^{(i)} \rangle \notin B(L)$ for all $i \geq 2$.

The directed analogue of Proposition ?? holds thanks to our directed connectivity assumption on $L$. In particular, lexicographic ordering should be replaced by the ordering induced by dot product with $\vec{v}$.

Finally, we conclude that the Ratio Limit Theorem ?? holds for directed bond animals, directed site animals, and directed bond trees, on any lattice of Section 3.1 (except for the Dead-End Lattice). This is because of the generalizations of Propositions ?? and ?? mentioned above, and by the fact that Theorem ?? does not distinguish between $(CA4)$ and $(DCA4)$, so the theorem and its proof are valid for directed clusters.

# 4 Proofs of Theorems

For real numbers $r$, we use the "floor" notation $\lfloor r \rfloor$ to denote the largest integer less than or equal to $r$.

**Proof of Theorem ??**: Before we begin, here is a very rough idea of the way the proof works. Since the pattern $P$ is bounded, there is a number $\rho > 0$ (depending on $P$) such that any cluster of size $n$ can be changed locally in $\lfloor \rho n \rfloor$ places to get $\lfloor \rho n \rfloor$ non-overlapping translates of $P$. But we will only choose $\lfloor \tau n \rfloor$ of these places to do this change, where $\tau$ is a fixed number between 0 and $\rho$. Now, let $G$ be a generic cluster of size $n$ that contains very few translates of $P$. If we choose $\lfloor \tau n \rfloor$ of the $\lfloor \rho n \rfloor$ possible locations mentioned above, then we get a cluster $H$ with $\lfloor \tau n \rfloor$ (or perhaps more) translates of $P$. Thus a single $G$ corresponds to $\binom{\rho n}{\tau n}$ different $H$'s. (Some of the $H$'s could be the same if we were unlucky enough to choose one of the places where a $P$ already existed, but this is not a big problem since there are not many $P$'s in $G$.) The quantity $\binom{\rho n}{\tau n}$ grows exponentially in $n$, and this would show that the number of $H$'s is exponentially larger than the number of $G$'s (which is what we want), except for the obvious problem that different $G$'s can give rise to the same $H$. How bad is this non-injectivity? Since $H$ has $\lfloor \tau n \rfloor$ translates of $P$ (or more, but not too many more), obtained by local changes, there are at most $K_1^{\tau n}$ $G$'s that give rise to $H$ (here and below, the $K_i$'s are constants). There are some other things that need to be taken into account too, including the size of $H$ (which need not be $n$, but is within $\pm \kappa \tau n$ of $n$), and the weight of $H$ (which is within a factor of $K_2^{\tau n}$ of $G$). Putting everything together, the weight of all $G$'s that get turned into $H$ is at most $K_3^{\tau n}$. So it comes down to a contest between $\binom{\rho n}{\tau n}$ and $K_3^{\tau n}$. Fortunately, however large $K_3$ is, we can choose $\tau$ small enough so that $\binom{\rho n}{\tau n}$ is exponentially larger than $K_3^{\tau n}$. And this is what we need to prove the theorem.

We now proceed with the proper proof. Since the set $D$ of Axiom $(CA4)$ can be enclosed in a finite box in $\mathbf{R}^d$, and since there is a finite upper bound on the number of sites in a unit hypercube of $\mathbf{R}^d$, it follows that there is a positive constant $\alpha$ such that for any $n$ and any $G \in C_n$, there exist at least $\lfloor \alpha n \rfloor$ sites $y_1, \ldots, y_{\lfloor \alpha n \rfloor} \in S(G)$ such that $(D + t(y_i)) \cap (D + t(y_j)) = \emptyset$ whenever $i \neq j$.

Consider a fixed $G \in C_n$, and fix the vectors $y_1, \ldots, y_{\lfloor \alpha n \rfloor}$ as described above. Next, given a number $\delta$ such that $0 < \delta < \alpha$, consider an arbitrary choice of $\lfloor \delta n \rfloor$ vectors from the set $\{y_1, \ldots, y_{\lfloor \alpha n \rfloor}\}$: call them $w_1, \ldots, w_{\lfloor \delta n \rfloor}$ (in some arbitrary order). Now define the sequence of clusters $G_0, \ldots, G_{\lfloor \delta n \rfloor}$ by

$$
\begin{aligned}
G_0 &= G \\
G_i &= T(G_{i-1}, w_i) \quad \text{for } i = 1, \ldots, \lfloor \delta n \rfloor.
\end{aligned} \tag{17}
$$

Let $H = G_{\lfloor \delta n \rfloor}$, and let $W$ be the (ordered) sequence $w_1, \ldots, w_{\lfloor \delta n \rfloor}$. By (??), the size of $H$ is between $n - \kappa \lfloor \delta n \rfloor$ and $n + \kappa \lfloor \delta n \rfloor$.

Consider an $\epsilon > 0$ (later, we shall set $\epsilon = \delta/2$). Consider the collection of all triples $(G, H, W)$ where $G$ is a cluster in $C_n^*$ that contains at most $\lfloor \epsilon n \rfloor$ translates of $P$, and where $H$ and $W$ are obtained from $G$ by the process described in the preceding paragraph. Considering the number of ways to choose the ordered sequence $W$, we obtain an immediate lower bound on the sum of $wt(G)$ over all such triples:

$$\sum_{(G,H,W)} wt(G) \geq \mathcal{G}_n[\leq \epsilon n, P] \frac{\lfloor \alpha n \rfloor!}{(\lfloor \alpha n \rfloor - \lfloor \delta n \rfloor)!} \ . \tag{18}$$

Now we shall derive a lower bound for this sum over triples. Notice that if $y \in S(G)$ and $G$ contains exactly $k$ translates of $P$, then it is possible that $T(G, y)$ contains more than $k + 1$ translates of $P$; however, any translate of $P$ that is in $T(G, y)$ but not in $G$ must overlap $D + t(y)$ (since everything outside of $D + t(y)$ is the same in both clusters). So let $q$ be the number of translates of $P_1 \cup P_2$ that intersect $D$; then we can be sure that $T(G, y)$ contains at most $k + q$ translates of $P$. Hence any $H$ in a triple $(G, H, W)$ contains at most $\epsilon n + q\lfloor \delta n \rfloor$ translates of $P$. Let $Z = 2^{(|S(D)|+|B(D)|)}$ be the number of subsets of $D$. Then for any cluster $K$ and any $w \in S(L)$, there are at most $Z$ clusters $G'$ such that $T(G', w) = K$. Also notice that if $K$ contains exactly $j$ translates of $P$, then there are at most $qj$ translates of $D$ that intersect one of these translates of $P_1 \cup P_2$. Hence there are at most $q|S(D)|j$ choices of $w$ in $S(L)$ for which $\{G' : T(G', w) = K\}$ is nonempty. Therefore, given $H$, there are at most $q|S(D)|(\lfloor \epsilon n \rfloor + q\lfloor \delta n \rfloor)Z$ ways to choose $w_{\lfloor \delta n \rfloor}$ and $G_{\lfloor \delta n \rfloor - 1}$ (recall (??)). If these are to be part of a valid triple, then $G_{\lfloor \delta n \rfloor - 1}$ contains at most $\lfloor \epsilon n \rfloor + q(\lfloor \delta n \rfloor - 1)$ translates of $P$, and so there are at most $q|S(D)|(\lfloor \epsilon n \rfloor + q(\lfloor \delta n \rfloor - 1))Z$ ways to choose $w_{\lfloor \delta n \rfloor - 1}$ and $G_{\lfloor \delta n \rfloor - 2}$. Therefore the number of triples $(G, H, W)$ in which any $H$ can occur is at most

$$\prod_{i=1}^{\lfloor \delta n \rfloor} q|S(D)|(\lfloor \epsilon n \rfloor + qi)Z \ ,$$

which in turn is less than

$$\frac{(\lfloor \epsilon n \rfloor + q\lfloor \delta n \rfloor)!}{(\lfloor \epsilon n \rfloor + q\lfloor \delta n \rfloor - \lfloor \delta n \rfloor)!} (q|S(D)|Z)^{\lfloor \delta n \rfloor} \ .$$

Together with (??), this implies that

$$\sum_{(G,H,W)} wt(G) \leq \sum_{(G,H,W)} \gamma^{\lfloor \delta n \rfloor} wt(H)$$

$$\leq \left( \sum_{j=n-\kappa\lfloor \delta n \rfloor}^{n+\kappa\lfloor \delta n \rfloor} \mathcal{G}_j \right) \frac{(\lfloor \epsilon n \rfloor + q\lfloor \delta n \rfloor)!}{(\lfloor \epsilon n \rfloor + q\lfloor \delta n \rfloor - \lfloor \delta n \rfloor)!} \Psi^{\lfloor \delta n \rfloor}, \tag{19}$$

where we have defined $\Psi = \gamma q |S(D)| Z$. Combining inequalities (??) and (??), taking $n^{th}$ roots, and letting $n \to \infty$, we find (using Stirling's formula) that

$$\limsup_{n \to \infty} (\mathcal{G}_n[\leq \epsilon n, P])^{1/n} \frac{\alpha^\alpha}{e^\delta (\alpha - \delta)^{\alpha - \delta}} \tag{20}$$
$$\leq \max\{\lambda^{1-\kappa\delta}, \lambda^{1+\kappa\delta}\} \frac{(\epsilon + q\delta)^{\epsilon + q\delta}}{e^\delta (\epsilon + q\delta - \delta)^{\epsilon + q\delta - \delta}} \Psi^\delta$$

(we need the "max" since $\lambda$ could be greater than or less than one). Recall that $\delta$ has not yet been specified, except that $0 < \delta < \alpha$. Write $\delta = \alpha t$, where $t$ is a number in $(0,1)$ that will be fully specified below. Also set $\epsilon = \delta/2$. Then we can rewrite (??) as

$$\frac{\limsup_{n \to \infty} (\mathcal{G}_n[\leq \epsilon n, P])^{1/n}}{\lambda} \left( \frac{1}{t^t (1-t)^{1-t}} \right)^\alpha \tag{21}$$
$$\leq \max\{\lambda^{-\kappa\delta}, \lambda^{+\kappa\delta}\} \Psi^\delta \left( \frac{(\frac{1}{2} + q)^{\frac{1}{2}+q}}{(\frac{1}{2} + q - 1)^{\frac{1}{2}+q-1}} \right)^\delta .$$

Now let

$$Q = \max\{\lambda^{-\kappa}, \lambda^{+\kappa}\} \Psi \frac{(\frac{1}{2} + q)^{\frac{1}{2}+q}}{(q - \frac{1}{2})^{q - \frac{1}{2}}} .$$

Then (??) becomes

$$\frac{\limsup_{n \to \infty} (\mathcal{G}_n[\leq \epsilon n, P])^{1/n}}{\lambda} \leq \left( t^t (1-t)^{1-t} Q^t \right)^\alpha , \tag{22}$$

which holds for every $t \in (0,1)$. Setting $t = 1/(1+Q)$ makes the right-hand side equal to $[Q/(Q+1)]^\alpha$, which is strictly less than 1. This proves the theorem.
□

Theorem ?? is a consequence of the following two propositions:

**Proposition 4.1** *Assume the Cluster Axioms (CA1) and (CA5), and assume that the conclusion of Theorem ?? holds. Then there is a positive constant $\Gamma$ such that*

$$\frac{\mathcal{G}_{n+2}}{\mathcal{G}_n} \geq \left( \frac{\mathcal{G}_{n+1}}{\mathcal{G}_n} \right)^2 - \frac{\Gamma}{n} \tag{23}$$

*for all sufficiently large $n$.*

**Proposition 4.2** *Let $\{a_n : n \geq 1\}$ be a sequence of positive numbers, and let $\phi_n = a_{n+1}/a_n$. Assume that there exist positive constants $\mu$ and $\Gamma$ such that*
*(a) $\lim_{n \to \infty} a_n^{1/n} = \lambda$,*
*(b) $\liminf_{n \to \infty} \phi_n > 0$, and*
*(c) $\phi_{n+1} \phi_n \geq (\phi_n)^2 - \Gamma/n$ for all sufficiently large $n$.*
*Then*

$$\lim_{n \to \infty} \phi_n = \lambda . \tag{24}$$

Notice that if $\Gamma$ were 0 on Proposition ??, then $\mathcal{G}_{n+1}/\mathcal{G}_n$ would be increasing in $n$. Thus Proposition ?? says roughly that if this sequence of ratios is almost monotone, and if $\mathcal{G}_n^{1/n}$ converges, then the sequence of ratios converges.

Proposition ?? is exactly Lemma 7.3.1 of Madras and Slade (1993) (except for the unimportant change of notation which changes the $n+2$ there to $n+1$ here; see the Remark following the proof of Lemma 7.3.1). We shall prove Proposition ?? below; it is similar to the proof of Theorem 7.3.2 in Madras and Slade (1993). Both of these results from Madras and Slade are in turn based on Kesten (1963). To deduce Theorem ??, take $a_n = \mathcal{G}_n$ in Proposition ??. Assumption (a) is axiom (CA3); assumption (b) is the assumption of $\Upsilon$ in Theorem ??; and assumption (c) is the conclusion of Proposition ??. The conclusion (??) is the conclusion of Theorem ??.

**Proof of Proposition ??:** For nonnegative integers $a$ and $b$, let $\mathcal{G}_n(a,b)$ denote the weighted sum of the set of clusters of size $n$ that contain exactly $a$ translates of $U$ and $b$ translates of $V$ (i.e., $\mathcal{G}_n(a,b)$ is the sum of $wt(G)$ over all $G \in C_n^*$ such that $|\tau_U(G)| = a$ and $|\tau_V(G)| = b$). Also let

$$\mathcal{G}_n(\geq a, \geq b) := \sum_{i \geq a, j \geq b} \mathcal{G}_n(i,j) \, .$$

In particular, $\mathcal{G}_n(\geq 0, \geq 0) = \mathcal{G}_n$.

First we note the identity

$$a\mathcal{G}_n(a,b) = \frac{b+1}{\theta}\mathcal{G}_{n+1}(a-1, b+1) \quad \text{for all } a \geq 1 \text{ and } b \geq 0. \tag{25}$$

To derive this identity, consider all pairs of clusters $(G, G')$ where $G \in C_n$, $|\tau_U(G)| = a$, $|\tau_V(G)| = b$, and $G' = \hat{G}_x$ for some $x \in \tau_U(G)$. Clearly, the left-hand side of the identity is the weighted sum of the set of all such pairs. But for such a pair we also have $G' \in C_{n+1}$, $|\tau_U(G')| = a-1$, $|\tau_V(G')| = b+1$, $G = \hat{G}'_x$ for $x \in \tau_V(G')$, and $wt(G) = wt(G')/\theta$ (by (CA5)). So the right-hand side of (??) also equals the weighted sum of all such pairs.

Using (??), we obtain

$$\begin{aligned}
\mathcal{G}_{n+1}(\geq 0, \geq 1) &= \sum_{i \geq 1, j \geq 0} \mathcal{G}_{n+1}(i-1, j+1) \\
&= \theta \sum_{i \geq 1, j \geq 0} \mathcal{G}_n(i,j) \frac{i}{j+1} \\
&= \theta \sum_{i \geq 0, j \geq 0} \mathcal{G}_n(i,j) \frac{i}{j+1} \tag{26}
\end{aligned}$$

and

$$\mathcal{G}_{n+2}(\geq 0, \geq 2) = \sum_{i \geq 2, j \geq 0} \mathcal{G}_{n+2}(i-2, j+2)$$

$$
\begin{align}
&= \theta^2 \sum_{i \geq 2, j \geq 0} \mathcal{G}_n(i,j) \frac{i(i-1)}{(j+1)(j+2)} \\
&= \theta^2 \sum_{i \geq 2, j \geq 0} \mathcal{G}_n(i,j) \frac{i(i-1)}{(j+1)(j+2)}. \tag{27}
\end{align}
$$

The Schwarz inequality implies

$$
\left( \sum \mathcal{G}_n(i,j) \frac{i}{j+1} \right)^2 \leq \left( \sum \mathcal{G}_n(i,j) \right) \left( \sum \mathcal{G}_n(i,j) \frac{i^2}{(j+1)^2} \right), \tag{28}
$$

where each sum is over $i \geq 1$ and $j \geq 0$. Inserting (??) into (??) gives

$$
\left[ \theta^{-1} \mathcal{G}_{n+1}(\geq 0, \geq 1) \right]^2 \leq \mathcal{G}_n \left( \sum_{i \geq 1, j \geq 0} \mathcal{G}_n(i,j) \frac{i^2}{(j+1)^2} \right), \tag{29}
$$

For $n \geq 1$, let

$$
\Xi_n = \frac{\mathcal{G}_{n+2}(\geq 0, \geq 2)}{\mathcal{G}_n} - \left( \frac{\mathcal{G}_{n+1}(\geq 0, \geq 1)}{\mathcal{G}_n} \right)^2
$$

and

$$
\Delta_n = \frac{\mathcal{G}_{n+2}}{\mathcal{G}_n} - \left( \frac{\mathcal{G}_{n+1}}{\mathcal{G}_n} \right)^2 - \Xi_n.
$$

Since $V$ is a proper pattern, Theorem ?? shows that the error term $\Delta_n$ decays to 0 exponentially rapidly as $n$ increases. So to prove the theorem it suffices to show $\Xi_n \geq -A/n$ for some constant $A$.

By (??) and (??),

$$
\begin{align}
\Xi_n &\geq \left( \sum_{i \geq 0, j \geq 0} \mathcal{G}_n(i,j) \frac{i(i-1)}{(j+1)(j+2)} - \sum_{i \geq 1, j \geq 0} \mathcal{G}_n(i,j) \frac{i^2}{(j+1)^2} \right) \frac{\theta^2}{\mathcal{G}_n} \\
&= \frac{\theta^2}{\mathcal{G}_n} \sum_{i \geq 0, j \geq 0} \mathcal{G}_n(i,j) \frac{(-i^2 - ij - i)}{(j+1)^2(j+2)} \tag{30}
\end{align}
$$

There exists a positive constant $K$ such that for every $n$, no cluster in $C_n$ contains more than $Kn$ translates of $U$ or of $V$. Hence the term $-i^2 - ij - i$ appearing in (??) is greater than $-3K^2 n^2$. Next, by Theorem ??, there exists $\epsilon > 0$ such that

$$
\limsup_{n \to \infty} \left( 1 - \frac{\mathcal{G}_n(\geq 0, \geq \epsilon n)}{\mathcal{G}_n} \right)^{1/n} < 1.
$$

Splitting the sum over $j$ in (??) into $\epsilon n \leq j \leq Kn$ and $0 \leq j < \epsilon n$, we obtain

$$
\Xi_n \geq \frac{-3K^2 n^2 \mathcal{G}_n(\geq 0, \geq \epsilon n)}{(\epsilon n)^3 \mathcal{G}_n} + (-3K^2 n^2) \left( 1 - \frac{\mathcal{G}_n(\geq 0, \geq \epsilon n)}{\mathcal{G}_n} \right).
$$

As $n \to \infty$, the first term in the right hand side is asymptotic to $-3K^2/n\epsilon^3$, and the second term decays to 0 exponentially. Thus the proposition is proven.
□

## 5 Discussion

This paper proves a Pattern Theorem and a Ratio Limit Theorem that hold for a wide variety of lattices and clusters, including bond animals, bond trees, and site animals, with different kinds of weights. The proofs were written to accommodate this generality, as well as to include other examples of lattices and clusters that other authors may need to consider. The most restrictive axiom is $(CA4)$, which does not seem to hold for self-avoiding walks, self-avoiding surfaces, or site trees. In addition, the Euclidean structure is important. For example, extending these results to graphs embedded in hyperbolic space (e.g., those of Swierczak and Guttmann (1996)) is not straightforward, partly because the notion of "translation" is no longer so simple. Nevertheless, it is reasonable to expect a version of the Pattern Theorem to hold for all of theses cases.

The pattern theorem suggests that there should be some kind of law of large numbers for pattern occurrence. That is, given a proper pattern $P$, does there exist a number $\alpha > 0$ such that "almost all" clusters of size $n$ contain between $(\alpha - \epsilon)n$ and $(\alpha + \epsilon)n$ translates of $P$? Or more simply, is the average number of patterns in a cluster of size $n$ asymptotically proportional to $\alpha n$? Nothing is known about these questions in general.

Another intriguing problem involves the universality of the critical exponent $\theta$, as described in the Introduction (recall Equation (??)). To fix ideas, consider bond animals on the square lattice $\mathbf{Z}^2$ and on the triangular lattice $Tri$. Using notation as in Corollaries ?? and ??, we believe that

$$|C^*_{BA,n}[Tri]| \sim K' n^{-\theta} \mu_{BA}[Tri]^n \quad \text{and} \quad |C^*_{BA,n}[\mathbf{Z}^2]| \sim K'' n^{-\theta} \mu_{BA}[\mathbf{Z}^2]^n \tag{31}$$

where $K'$ and $K''$ are positive constants, and $\theta$ has the same value in both expressions. It seems very hard to prove the relations of (??). An easier task might be to prove the following consequence of (??):

$$\frac{|C^*_{BA,n}[\mathbf{Z}^2]|}{|C^*_{BA,n}[Tri]|} \sim K \rho^n, \tag{32}$$

where $K$ is a constant and $\rho = \mu_{BA}[\mathbf{Z}^2]/\mu_{BA}[Tri] < 1$. As in Corollary ??, think of bond animals of $\mathbf{Z}^2$ as bond animals of $Tri$ that contain no diagonal bond (recall Figure ??). Let $P$ the pattern $(b, \emptyset)$, where $P$ is a diagonal bond. Then Equation (??) says that the fraction of bond animals on $Tri$ that contain <u>no</u> translates of $P$ decays purely exponentially, with no multiplicative power law term. A proof of this assertion would be very strong support for the universality of $\theta$, even in the absence of a rigorous proof that $\theta$ exists. And it is conceivable that information about the simpler expression in (??) may be more accessible than information about (??).

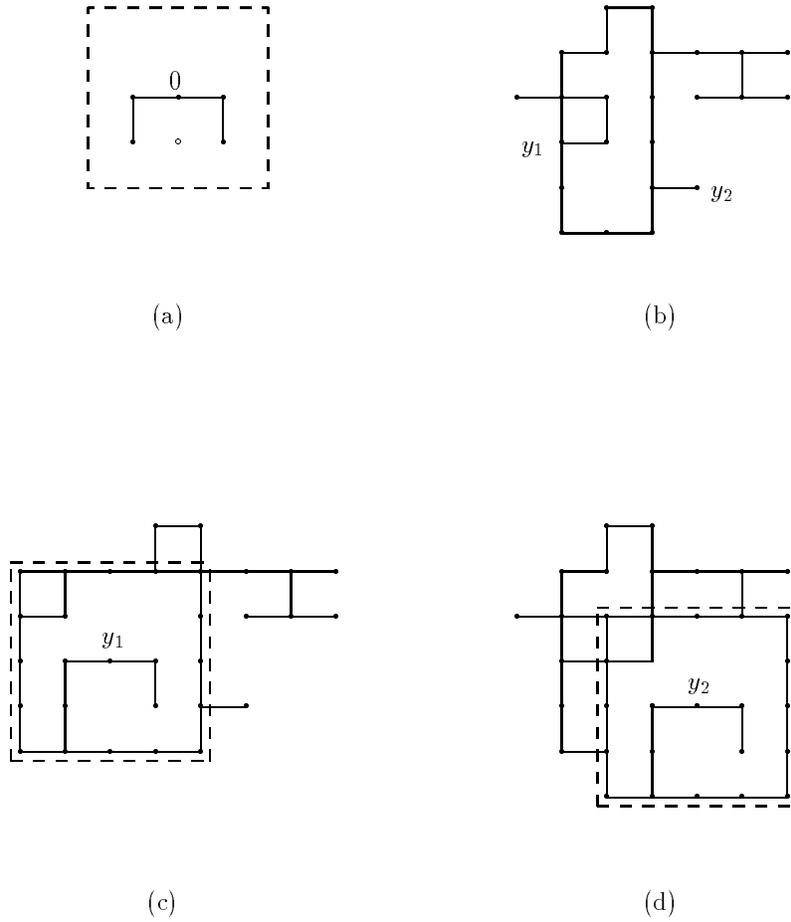

Figure 1: Cluster Axiom ($CA4$) in $\mathbf{Z}^2$: (a) A pattern $P = (P_1, P_2)$: $P_2$ is the single site marked by the circle; $P_1$ is the solid lines and dots (four bonds and five sites); the dashed square is the boundary of $D$; 0 is the origin. (b) A bond animal $G$ with two sites labelled. (c) One possibility for $T(G, y_1)$; the dashed square is included to surround $D + t = D + y_1$. (d) One possibility for $T(G, y_2)$. See also Figure ??.

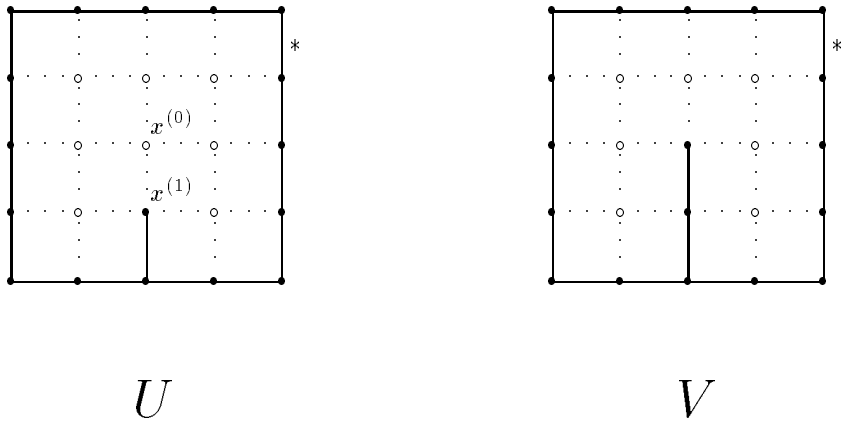

Figure 2: Examples of patterns $U$ and $V$ for Axiom ($CA5$) in $\mathbf{Z}^2$, as well as for their $(1,1)$-directed versions. See Proposition **??** in Section 3.5 for the meaning of $x^{(0)}$ and $x^{(1)}$. Circles and dotted lines denote sites and bonds of $U_2$ and $V_2$; solid dots and lines denote $U_1$ and $V_1$. To obtain $U$ and $V$ for bond trees (including $(1,1)$-directed trees), move the bond labeled '*' from $U_1$ to $U_2$ and from $V_1$ to $V_2$.

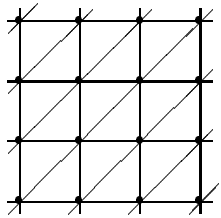
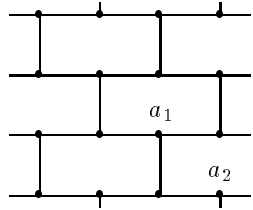

Triangular            Hexagonal

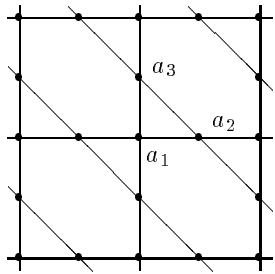
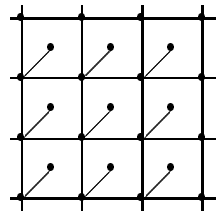

Kagome            Dead-End

Figure 3: Some examples of lattices from Section 3.1.

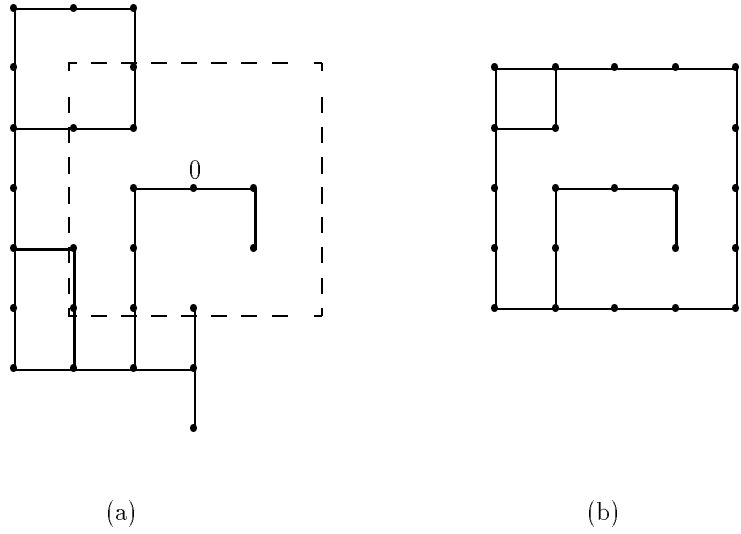

(a)                          (b)

Figure 4: Proof of Proposition ??: (a) A bond animal $H$ that contains the pattern $P$ of Figure ??(a). The dashed line indicates $\partial D$, the boundary of $D$. (b) The resulting $\tilde{H} = (H \cap D) \cup \partial D$. The animals $T(G, y_1)$ and $T(G, y_2)$ in Figure ??(c,d) were obtained using this $\tilde{H}$.

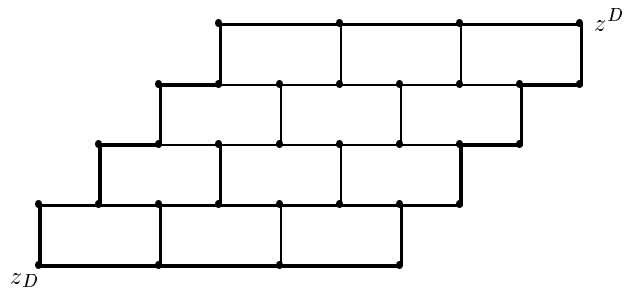

Figure 5: From the proof of Proposition ??: A possible set $D$ for $(1,1)$-directed clusters on the hexagonal lattice. The thicker lines denote $\partial D$.